# A *virtual* substitution of Brouwer's choice sequence

Dipl. Math.(FH) Klaus Lange*

*Cinterion Wireless Modules GmbH, Technology, Product Development, Integration Test
Berlin, Germany
Klaus.Lange@cinterion.com

______________________________________________________________________

Abstract: Step by step a substitution of the well known Brouwer's choice sequence ([1] to [3]) will be constructed. It begins with an establishing of quasi alternating prime number series followed by a construction of a virtual sequence in sense of the virtual set definition [5]. The last step gives reasons for why this virtual sequence substitutes the choice sequence created by L. E. J. Brouwer.

______________________________________________________________________

## I. Establishing quasi alternating prime number series

Firstly we are going to establish quasi alternating prime number series. Quasi alternating means that we use different forms of prime numbers and change the leading sign during summation.

Be e an even natural number and o the odd remain from 1 up to e – 1 we write for the prime numbers p the formula

$$p = en + o \qquad (1)$$

with $n \in \mathbb{N}$.

We search for all e in (1) that fits for two preconditions:

   a) there is only an even number of forms which fulfil e|o = 1 for given e
   b) if we use an even number of pairs with prime number forms (with different leading sign in every pair) then for the infinite quasi alternating series

$$S = \sum_{i \to \infty} m_i \qquad (2)$$

- with $m_i$ as prime number forms with different leading signs - infinite many pairs of part sums of S like

$S_1 = m_1, S_2 = m_1+m_2, S_3 = m_1+m_2+m_3, \ldots, S_j, S_{j+1}, \ldots$ exist, that

$S_j > 0 \land S_{j+1} < 0$.

Some examples:

<u>e = 4:</u>

o = 1 and o+2 = 3 fit

$-p = -(4n + 1) \land q = 4n + 3$

leads to

$S = 3 - 5 + 7 + 11 - 13 - 17 + 19 + 23 - 29 + 31 - 37 + \ldots$

$S_1 = 3, S_2 = -2, S_3 = 5, S_4 = 16, S_5 = 3, S_6 = -14, S_7 = 5, S_8 = 28, S_9 = -1, \ldots$

<u>e = 6:</u>

o = 1 and o = 5 fit

$-p = -(4n + 1) \land q = 4n + 5$

because for o = 3 it is 6|o = 3 and we have only two forms.

$S = 5 - 7 + 11 - 13 + 17 - 19 + 23 + 29 - 31 - 37 + 41 - 43 + \ldots$

$S_1 = 5, S_2 = -2, S_3 = 9, S_4 = -4, S_5 = 13, S_6 = -6, S_7 = 17, S_8 = 46, S_9 = 15, S_{10} = -22, \ldots$

<u>e = 8:</u>

o = 1, o = 3, o = 5, o = 7 fit, because all of them fulfil 8|o = 1.

we have a lot of possible forms, so we can change the leading sign:

$$-p = -(8n + 1) \land -q = -(8n + 3) \land u = 8n + 5 \land v = 8n + 7$$

or

$$-p = -(8n + 1) \land q = 8n + 3 \land -u = -(8n + 5) \land v = 8n + 7$$

or

or

$$-p = -(8n + 1) \wedge q = 8n + 3 \wedge u = 8n +5 \wedge -v = -(8n +7)$$

or

$$p = 8n + 1 \wedge -q = -(8n + 3) \wedge u = 8n +5 \wedge -v = -(8n +7)$$

or

$$p = 8n + 1 \wedge -q = -(8n +3) \wedge -u = -(8n + 5) \wedge v = 8n +7$$

or

$$p = 8n + 1 \wedge q = 8n + 3 \wedge -u = -(8n +5) \wedge -v = -(8n +7)$$

with six different series and part sums. So we have more alternatives to build the series which satisfy the precondition b).

e = 10:

With o=1, o=3, o=7 and o=9 we have four prime number forms and we have more than one alternative to build the series under given preconditions. In the same way as in example e=8.

e = 12:

In this case we have only four useable prime number forms because for o=3 and o=9 it is o|12= 3. The four prime number forms $12n + o$ with $o = \{1, 5, 7, 11\}$ have a lot alternatives to change the leading sign as in e = 8.

e = 14:

Here we have six useable prime number forms, only for o=7 it is e|o > 1. All possible combinations building alternatives for quasi prime number series will shown in a table. We see 20 alternatives setting the leading signs of p for $^j m_i$:

| 14n + o | | Number j of alternatives building quasi prime number series a for $^j m_i$ | | | | | | | | | | | | | | | | | | | |
|---|---|---|---|---|---|---|---|---|---|---|---|---|---|---|---|---|---|---|---|---|---|
| o | | 1 | 2 | 3 | 4 | 5 | 6 | 7 | 8 | 9 | 10 | 11 | 12 | 13 | 14 | 15 | 16 | 17 | 18 | 19 | 20 |
| 1 | | + | - | - | - | + | + | + | + | + | + | - | - | - | - | - | - | + | + | + | - |
| 3 | | + | + | + | + | - | - | - | + | + | + | - | - | - | + | + | + | - | - | - | - |
| 5 | | + | + | + | + | + | + | + | - | - | - | + | + | + | - | - | - | - | - | - | - |
| 9 | | - | + | - | - | + | - | - | + | - | - | + | + | - | + | + | - | + | + | - | + |
| 11 | | - | - | + | - | - | + | - | - | + | - | + | - | + | + | - | + | + | - | + | + |
| 13 | | - | - | - | + | - | - | + | - | - | + | - | + | + | - | + | + | - | + | + | + |

Table 1: Alternatives building quasi alternating prime number series for e = 14

After this examples we generalise the building of alternatives for quasi alternating prime number series.

Using (1) for e >> 0 we can try to sort the alternatives using $m_i$ counting them in formations $^j m_i$ in this way:

Precondition for all entries which are used in (3a) and (3b) is that e and o in (1) satisfy e|o=1. The number of entries for –p and p must be the same.

$$^1m_i = \begin{cases} -p, \text{ if } p \in \{en + o\} \text{ with } o = 1 \\ -p, \text{ if } p \in \{en + o\} \text{ with } o = 3 \\ -p, \text{ if } p \in \{en + o\} \text{ with } o = 5 \\ \ldots \\ -p, \text{ if } p \in \{en + o\} \text{ with } o = e/2 - 5 \\ -p, \text{ if } p \in \{en + o\} \text{ with } o = e/2 - 3 \\ -p, \text{ if } p \in \{en + o\} \text{ with } o = e/2 - 1 \\ \\ p, \text{ if } p \in \{en + o\} \text{ with } o = e/2 + 1 \\ p, \text{ if } p \in \{en + o\} \text{ with } o = e/2 + 3 \\ p, \text{ if } p \in \{en + o\} \text{ with } o = e/2 + 5 \\ \ldots \\ p, \text{ if } p \in \{en + o\} \text{ with } o = e - 5 \\ p, \text{ if } p \in \{en + o\} \text{ with } o = e - 3 \\ p, \text{ if } p \in \{en + o\} \text{ with } o = e - 1 \end{cases} \quad (3a)$$

to

$$^Zm_i = \begin{cases} p, \text{ if } p \in \{en + o\} \text{ with } o = 1 \\ p, \text{ if } p \in \{en + o\} \text{ with } o = 3 \\ p, \text{ if } p \in \{en + o\} \text{ with } o = 5 \\ \ldots \\ p, \text{ if } p \in \{en + o\} \text{ with } o = e/2 - 5 \\ p, \text{ if } p \in \{en + o\} \text{ with } o = e/2 - 3 \\ p, \text{ if } p \in \{en + o\} \text{ with } o = e/2 - 1 \\ \\ -p, \text{ if } p \in \{en + o\} \text{ with } o = e/2 + 1 \\ -p, \text{ if } p \in \{en + o\} \text{ with } o = e/2 + 3 \\ -p, \text{ if } p \in \{en + o\} \text{ with } o = e/2 + 5 \\ \ldots \\ -p, \text{ if } p \in \{en + o\} \text{ with } o = e - 5 \\ -p, \text{ if } p \in \{en + o\} \text{ with } o = e - 3 \\ -p, \text{ if } p \in \{en + o\} \text{ with } o = e - 1 \end{cases} \quad (3b)$$

Notes:
(i) Z is the last alternative to build a quasi alternating prime number series for given e.
Z is determinated by the given e and the number of terms for p and –p.
That depends on the precondition e|o = 1.
The complete entries for all odd remains in (3a) and (3b) o = {1, 3, 5, …, e-5, e-3, e-1} we only have if $e = 2^a$ with a > 1 and a ∈ IN. **Appendix A** shows an example of $e = 16 = 2^4$ with 8 odd remains o = {1, 3, 5, …, 11, 13, 15}. This few remains cause 70 different alternatives building quasi alternating prime number series!
(ii) Why we can use the precondition b) for crossing the zero line infinite many times? The reason is that all prime number forms which fit precondition a) satisfy the theorem of P. G. L. Dirichlet that means that all prime number forms en + o with e|o = 1 exist infinite many times and statistically are in equilibrium. Effects known as 'Chebyshev's bias' for example the form pair (4n + 1, 4n + 3) [4] will be balanced out by an extra change of the leading sign during summation. Later on we will analyse this topic and give a solution.

## II. Searching for zeros

If we talk about quasi alternating series the question to search for zeros is near. That means if we have the series

$$S = \sum_{i=1 \to \infty} m_i$$

then we ask for sums

$$S_n = \sum_{i=1 \to n} m_i = 0. \tag{4}$$

This is equal to

$$S_{n-1} = \sum_{i=1 \to n-1} m_i = -m_n \tag{4a}$$

$$S_{n-2} = \sum_{i=1 \to n-2} m_i = -(m_{n-1} + m_n)$$

$$S_{n-3} = \sum_{i=1 \to n-3} m_i = -(m_{n-2} + m_{n-1} + m_n)$$

…..

$$S_1 = m_1 = -\sum_{i=2 \to n} m_i$$

After reaching the zero with $m_n$ the series could restart with $m_{n+1}$ and run up to the next zero and so on.

But do every of these series S have a part sum with value zero?

This question challenges all aspects of the analysed topic. The try to find an answer leads to a new approach and gives another point of view to the foundations of formalised modern mathematics. To give a first answer we look ones more to a known example.

Using (1) as example with e = 6 and o = {1; 5} and (3a)

$$m_i = \begin{cases} -p, & \text{if } p \in \{6n + 1\} \\ p, & \text{if } p \in \{6n + 5\} \end{cases}$$

we find

S = 5 –7 + 11 – 13 + 17 – 19 + 23 + 29 – 31 – 37 + 41 – 43 + 47 + 53 + 59 – 61 – 67 + …

and for this series the part sums

$S_{78}$ = 5 – 7 + 11 – 13 + …. – 367 – 373 – 379 + 383 + 389 – 397 + 401 – 409 = 0

$S_{84}$ = 5 – 7 + 11 – 13 + …. – 409 + 419 – 421 + 431 – 433 – 439 + 443 = 0

$S_{528}$ = 5 – 7 + 11 – 13 + …. – 3793 + 3797 + 3803 + 3821 = 0

One of these part sums with value zero is enough.

But for some prime number formats it is not easy to find a zero part sum.

For the first prime number format with e = 4 the search for zeros is difficult.

Using (3a) for e = 4

$$m_i = \begin{cases} -p, & \text{if } p \in \{4n + 1\} \\ p, & \text{if } p \in \{4n + 3\} \end{cases}$$

then we have

S = 3 – 5 + 7 + 11 – 13 – 17 + 19 + 23 – 29 + …

During our search for zeros we find a lot of changes for the leading sign of part sums but no zero.

Some examples will illustrate this behaviour:

$S_{195}$ = 821 > 0 > $S_{196}$ = -380

$S_{1216}$ = -9844 < 0 < $S_{1217}$ = 27

$S_{5633}$ = 40976 > 0 > $S_{5634}$ = -14553

$S_{9908}$ = -63855 < 0 < $S_{9909}$ = 39848

$S_{13675}$ = 76146 > 0 > $S_{13676}$ = -71915

With brute force computing it is possible to find a zero for large $S_i$ but we have no way to be sure finding a zero or not.

## III. Analysing the quasi alternating prime number series

Firstly the main problem is to find a zero for those given prime number series by an analytical prediction. This problem is determined by the prime number sequence itself.

Predict a zero for this series means that we try to predict the next prime number $p_{n+1}$ only by knowing the prime numbers $p_1$ to $p_n$ and without factorisating all numbers greater $p_n$. And this is not possible.

In (4a) we see that the problem is much stronger than that.

$$\sum_{i=1 \to n-1} m_i = - m_n$$

On the left side of the equation we see a sum of all prime numbers starting with 3 or 5 up to a given one. And on the right side we see only one prime number (with changed leading sign), which is the direct next greater of the given one.

Using the last example with $p = 4n \pm 1$ and the series starting with 3 we find

$S_3 = 3 - 5 + 7 = 5$

The part sum of the first three prime numbers starting with 3 has the value 5 and this is a prime number, but not 11 (the next after 7).

$S_{11} = 3 - 5 + 7 + 11 - 13 - 17 + 19 + 23 - 29 + 31 - 37 = -7$

This sum has as result the prime number 7 (with negative leading sign). But this is smaller than the next prime number after 37.

$S_{195} = 3 - 5 + 7 + 11 - 13 - 17 + \ldots + 1171 - 1181 + 1187 - 1193 = 821$

The next prime number after 1193 is 1201 and this is larger than the sum value 821. But 821 is a prime number too.

A larger part sum endens our examples

$S_{319} = 3 - 5 + 7 + 11 - 13 - 17 + \ldots - 2069 - 2081 + 2083 + 2087 - 2089 + 2099 + 2111 - 2113 - 2129 = 7457$

This is a prime number but it is greater than 2131, the direct follower of 2129.

A table summarizes the results:

| $m_{n-1}$ | $m_n$ | $S_{n-1}$ |
|---|---|---|
| 7 | 11 | 5 |
| -37 | -41 | -7 |
| -1193 | -1201 | 821 |
| -2129 | 2131 | 7457 |

Table 2: Relations between leading signs for $m_{n-1}$ and $S_{n-1}$

Finding zeros it is not only important that $m_n$ and $S_{n-1}$ have the same absolute value. Additional the prime number $m_n$ needs to have the correct prime number form, so that the leading signs of $S_{n-1}$ and $m_n$ are different.

In other words: We try to predict the prime number form of $m_n$ only by using all prime numbers smaller $m_n$. And that is impossible, either.

Secondly it could be possible that a boundary $\beta$ exists for which all values $S_i$ are greater or smaller then zero if $i \in \mathbb{N}$ is large enough in this way:

$$S_i < S_\beta < 0 \quad \vee \quad 0 < S_\beta < S_i \quad \text{with } i > \beta \qquad (5)$$

The reason for this we talked about in Note (ii). If we find for some prime number forms (1) the same effect called 'Chebyshev's bias' like for $p = 4n \pm 1$ we have to rule out that those properties for quasi alternating series could lead to a boundary $\beta$. The existence of a $\beta$ makes it in principle possible to calculate all part sums $S_1$ up to $S_{\beta-1}$ and search for zeros, and from $S_\beta$ and above in (5) we know that no more zeros exist. If those $\beta$ is calculated then in finite calculations we can decide about the existence of zeros.

How can the value of $\beta$ be calculable? That is an open question, but the possibility for some special forms (1) can not be ruled out. It is not enough that by heuristically argumentations some approximations can be made and also show that a $\beta$ must exist, without a calculation of the explicit value of $\beta$ is done. In this case we do not know how many calculation steps we have to made when searching for zeros.

The problem here is to decide between a membership of a subset for all series with a zero and without a zero. It is clear. If we have a prime series with a $\beta$ it could be possible that before reaching this $\beta$ we find a zero. How long or how many steps we must do while searching for a zero if we do not have a calculated value for $\beta$? That is another problem.

Back to our discussed main problem: To construct a virtual sequence without decidable quasi alternating prime series we need a mechanism that makes sure that infinity changes of leading signs in part sums happen.

Trying this we have to expand the way how to build alternatives of quasi alternating prime number series as seen in (3a) and (3b). The only reason for those constructions is to be sure that those prime number series without a $\beta$ in order of (5) could exist and that we could use those additional alternatives to construct a virtual sequence.

# IV. Construction of a virtual sequence

First of all we remind of the definition of a virtual subset given in [5].

<u>Definition Virtual subset</u>

A virtual subset $V \subset W$, $W \neq \{\}$ is given by the following criterions

a) $W \neq V$
b) $W = U_1 \cup U_2 \cup \ldots \cup U_n \cup V$ with $U_i \cap U_j = \{\}$ and $U_i \cap V = \{\}$.
c) In the sense of Gödel theorem [6] it is not decidable that either one element of V exists or that $V = \{\}$.

If $\mathbf{E} = \{$ e from prime number form $p = en + o\} = \{4; 6; 8; 10; \ldots\}$ then we have subsets of $\mathbf{E}$ for those e in (1) with all alternatives of constructible quasi alternating prime number series which have a zero. We call this subset $\mathbf{E_0}$.

The example of a member of $\mathbf{E_0}$ is the number $e = ,6$. For both alternatives to construct the quasi alternating prime number series in the sense of (3a) and (3b)

$$^1m_i = \begin{cases} -p, \text{ if } p \in \{6n + 1\} \\ p, \text{ if } p \in \{6n + 5\} \end{cases}$$

or

$$^2m_i = \begin{cases} p, \text{ if } p \in \{6n + 1\} \\ -p, \text{ if } p \in \{6n + 5\} \end{cases}$$

we find a zero.

Those e which have not only quasi alternating prime number series with zeros but those with a provable β – for all alternatives without zeros - are members of another subset of $\mathbf{E}$ and we will call it $\mathbf{E_β}$.

One special example of those members of $\mathbf{E_β}$ is a constructed prime number series with a β described in (5) after that we know that no more changes of the leading signs happen for the part sums $S_i$.

For

$$^1m_i = \begin{cases} -p, \text{ if } p \in \{4n + 1\} \\ p, \text{ if } p \in \{4n + 3\} \end{cases}$$

we observe for very large randomly picked x there is, with high possibility, a $p \leq x$ that has the form $p = 4n + 3$ instead of form $4n + 1$. So the positive prime numbers in the series win the battle and it exists a $\beta$. On the other hand we cannot find a zero for all $|m_i| < 10^n = 10^4$. Brute force computing is needed to verify more $|m_i|$ larger than $10^4$. So $e = 4$ could be an example of a membership of $\mathbf{E_\beta}$.

It is clear that the same would be correct for

$$^2m_i = \begin{cases} p, & \text{if } p \in \{4n + 1\} \\ -p, & \text{if } p \in \{4n + 3\} \end{cases}$$

Here we find an interesting point. What happens if we have a construction of a quasi alternating prime number series for a given e which have infinite many changes in the leading signs of the part sums $S_i$, that means there is no $\beta$, without finding a zero? We can not say that e is a member of $\mathbf{E_0}$ or a member of $\mathbf{E_\beta}$.

We have

$$\mathbf{E_0} \cap \mathbf{E_\beta} = \{\} \qquad (6)$$

But it is

$$\mathbf{E} \geq \mathbf{E_0} \cup \mathbf{E_\beta} \qquad (7)$$

We need a third kind of a subset. With this third subset we complete the last equation to

$$\mathbf{E} = \mathbf{E_0} \cup \mathbf{E_\beta} \cup \mathbf{E_\infty} \qquad (8)$$

$\mathbf{E_\infty}$ is a virtual subset, because we do not know if a subset $\mathbf{E_\infty}$ exist. Firstly it is possible that all alternatives $^xm_i$ for given e in a quasi alternating prime number series have a boundary $\beta$. And secondly, if this is not the case, those alternatives could have a zero for a very large index i.

The foundation of the existence of $\mathbf{E_\infty}$ is being sure that quasi alternating prime number series with infinite many changes of the leading sign in (2) exists. Those special series we can construct out of those of $\mathbf{E_\beta}$ in a simply way. One example of how to construct a member of $\mathbf{E_\infty}$ will be this:

Be $e \in \mathbf{E_\beta}$ with the calculated boundary $\beta$ for the special construction alternative like in (2)

$$^xS = \sum_{i=1 \to \infty} {}^xm_i$$

then it exists a special construction alternative with

$$^{y}S = \sum_{i=1\to\infty} {}^{y}m_i = \sum_{i=1\to\infty} - {}^{x}m_i$$

and the leading signs of all of $m_i$ are changed from $^{x}S$ to $^{y}S$.

Knowing this we will construct a quasi alternating prime number series with infinite many changes in leading signs for the part sums

$$^{\infty}S = \sum_{i=1\to\beta} {}^{x}m_i + \sum_{i=\beta+1\to 2\beta} {}^{y}m_i + \sum_{i=2\beta+1\to 3\beta} {}^{x}m_i + \sum_{i=3\beta+1\to 4\beta} {}^{y}m_i + \ldots \quad (9)$$

In (9) the variable $d\beta$ means $d\beta := \{\beta_1, \beta_2, \beta_3, \ldots\}$ with $d = \{1, 2, 3, \ldots\}$ because the explicit values of boundaries $\beta$ - ignoring they multitudes - for the several subsumes could be different.

The new constructed series in (9) now could have a zero because the new $^{\infty}S_n$ changes infinite many times their leading signs, so the possibility of finding a zero is given. But the only way to test this is to calculate these part sums $^{\infty}S_n$ step by step. Because $^{\infty}S$ is an additional series and not a substitute for the genuine S in (2) with calculated $\beta$, we have (8) confirmed. If we cannot find a zero during brute force computing for the series described in (9) the e of given basic primes forms (1) is a candidate for $\mathbf{E}_\infty$. But we cannot be sure of this membership which is a property of a virtual subset of $\mathbf{E}$. The confirmation of (8) leads to expand (6) by adding

$$\mathbf{E}_0 \cap \mathbf{E}_\infty = \{\} \quad (6a)$$

and

$$\mathbf{E}_\beta \cap \mathbf{E}_\infty = \{\} \quad (6b)$$

With (6) to (9) we have established a virtual subset as it is defined in [5].

Now it is easy to generate a virtual sequence.

Be $^{\infty}e_1$ the smallest element of $\mathbf{E}_\infty$, $^{\infty}e_2$ the second smallest of $\mathbf{E}_\infty$ and so on, then we can write

$$^{\infty}e_1 \; ^{\infty}e_2 \; ^{\infty}e_3 \ldots \text{ is a virtual sequence noted as } <^{\infty}e_i>i \quad (10)$$

That gives us the key to make the last step to substitute the choice sequence of Luitzen Egbertus Jan Brouwer.

Note:
(iii)   The more the number of alternatives of a quasi alternating series rises, we have

more possibilities to build different $^\infty S$ like in (9). If a given alternative $^x m_i$ leads to a boundary β with, for example, $^x S_i < {}^x S_\beta < 0$, it is not necessary that we use $^y m_i$ with the special property $^y m_i = - {}^x m_i$ but every other alternative $^z m_i$ whose boundary β satisfy $0 > {}^z S_\beta > {}^z S_i \ \forall i > \beta$. Once more β means $\beta := \{\ldots, \beta_x, \beta_y, \beta_z, \ldots\}$. Tables in the **Appendix A** give an impression of what this circumstance means for the risen number of possibilities building those kinds of $^\infty S$ depending on e in (1).

# V. Substitution of Brouwers choice sequence

L. E. J. Brouwer tries to find an example of breaking the 'Tertium non datur' in formalistic mathematics using his so called choice sequence [1] . Doing this he argues with parameters by free choice in this way:

Be $a_0 \ a_1 \ a_2 \ a_3 \ \ldots$ a choice sequence noted as $\langle a_i \rangle i$

then we look at the statement

$$\exists i (a_i = 0) \ \lor \ \forall i \ (a_i \neq 0).$$

There is no reason that this statement is true or false because of the free choice for all $a_i$ which is made in future. We have open parameters by free will without any law limitations. Those so called 'choice sequences' are not really sequences in a mathematical defined sense. In fact the term 'choice sequence' is more like an expression. This is shown in [2] and [3] with the result that every formula without choice parameters is equivalent with a formula that has only lawful function parameters.

That means: The problem here is not the question of 'Tertium non datur', but the question of lawful mathematical objects allowing a systematic analysing.

To make a step forward to the question of 'Tertium non datur' we need those lawful sequences without the possibility to predict if a given number is an element of this sequence or not.

Be $^4 E$ the set of even natural numbers e with $e > 2$.

In the last chapter we saw that it is possible to construct for every $e \in {}^4 E$ a series $^\infty S$ as it is given in (9). Only calculations can decide if we have a zero for this constructed series or not.
Now we can define a 'candidate set' $C \subset {}^4 E$

$C := \{c | \exists \ ^\infty S$ of c and up to this moment brute force computing cannot find a zero for $^\infty S\}$.

Using (10) we look to the statement

$$\exists i (^\infty e_i = c) \ \lor \ \forall i \ (^\infty e_i \neq c).$$

Threw brute force computing we only could rule out that c is not a member of $<^\infty e_i>i$ if a zero is found. And only the future will tell this after adding more calculations for $^\infty S_n$. For this prime number series no recursive formula can be made. A positive statement of membership in $<^\infty e_i>i$ needs infinite many calculations.

The difference about Brouwers choice sequence $<a_i>i$ and the virtual sequence $<^\infty e_i>i$ is that for the virtual one we do not need a free will, a free choice. The virtual sequence is determinate by well defined calculation steps. Both sequences have in common that only the future calculations can make a decision if c is not a member. But it never can be made a positive decision if c is really a member of $<^\infty e_i>i$. On one hand that circumstance challenges the 'Tertium non datur' in Brouwers sense and on the other hand the properties of [2] and [3] having only lawful determined parameters are now satisfied.

Special note: Conclusion and outlook

(iv)  Using constructed virtual subsets or virtual sequences like (10) the objection comes that *regardless whether our ability to decide* an e exists either in this or in another of the given subsets from **E**. In my opinion those statements make implicit preconditions. In other words: 'Tertium non datur' has become an implicit axiom. Virtual subsets and a virtual sequence as a substitution of Brouwers choice sequence remind us to analyse those implicit axioms or preconditions. The correctness of those axioms must be questioned. Doing this we could look at the results of Kurt Gödel gave 1938 [7]. For us, the main interest in this work is the result that if only defined constructed sets are allowed and not other pure formalistic mathematical objects then the axiom of choice is not an axiom but a provable and true theorem. That corresponds to aspects of Brouwer' s works [1], that only those mathematical objects exists or are correctly and consistently useable, which are given in a constructive way. In my next paper I will highlight this topic. For now in **Appendix B** we have a short analogy to quantum physics in understanding this new approach for an old discussion in mathematics.

# Appendix A

All possible combinations building alternatives of quasi prime number series for e = 16 in (1) are shown in the following table, which must be separated in parts to show all 70 alternatives setting the leading signs of p for $^jm_i$. This illustrates what happens alone for special case e=$2^a$ with a > 3 to find more different $^\infty S$ given in (9) having been discussed in note (iii).

| 16n + o | Number of alternatives for $^jm_i$ | | | | | | | | | | | | | | | | |
|---|---|---|---|---|---|---|---|---|---|---|---|---|---|---|---|---|---|
| o | 1 | 2 | 3 | 4 | 5 | 6 | 7 | 8 | 9 | 10 | 11 | 12 | 13 | 14 | 15 | 16 | 17 |
| 1 | + | + | + | + | + | + | + | + | + | + | + | + | + | + | + | + | + |
| 3 | + | + | + | + | + | + | + | + | + | + | + | + | + | + | + | - | - |
| 5 | + | + | + | + | + | - | - | - | - | - | - | - | - | - | - | + | + |
| 7 | + | - | - | - | - | + | + | + | + | - | - | - | - | - | - | + | + |
| 9 | - | + | - | - | - | + | - | - | - | + | + | + | - | - | - | + | - |
| 11 | - | - | + | - | - | - | + | - | - | + | - | - | + | + | - | - | + |
| 13 | - | - | - | + | - | - | - | + | - | - | + | - | + | - | + | - | - |
| 15 | - | - | - | - | + | - | - | - | + | - | - | + | - | + | + | - | - |

Table 3a: Alternatives No. 1 to 17 for e = 16

| 16n + o | Number of alternatives for $^jm_i$ | | | | | | | | | | | | | | | | | |
|---|---|---|---|---|---|---|---|---|---|---|---|---|---|---|---|---|---|---|
| o | 18 | 19 | 20 | 21 | 22 | 23 | 24 | 25 | 26 | 27 | 28 | 29 | 30 | 31 | 32 | 33 | 34 | 35 |
| 1 | + | + | + | + | + | + | + | + | + | + | + | + | + | + | + | + | + | + |
| 3 | - | - | - | - | - | - | - | - | - | - | - | - | - | - | - | - | - | - |
| 5 | + | + | + | + | + | + | + | + | - | - | - | - | - | - | - | - | - | - |
| 7 | + | + | - | - | - | - | - | - | + | + | + | + | + | + | - | - | - | - |
| 9 | - | - | + | + | + | - | - | - | + | + | + | - | - | - | + | + | + | - |
| 11 | - | - | + | - | - | + | + | - | + | - | - | + | + | - | + | + | - | + |
| 13 | + | - | - | + | - | + | - | + | - | + | - | + | - | + | + | - | + | + |
| 15 | - | + | - | - | + | - | + | + | - | - | + | - | + | + | - | + | + | + |

Table 3b: Alternatives No. 18 to 35 for e = 16

| 16n + o | Number of alternatives for $^jm_i$ | | | | | | | | | | | | | | | | | |
|---|---|---|---|---|---|---|---|---|---|---|---|---|---|---|---|---|---|---|
| o | 36 | 37 | 38 | 39 | 40 | 41 | 42 | 43 | 44 | 45 | 46 | 47 | 48 | 49 | 50 | 51 | 52 | 53 |
| 1 | - | - | - | - | - | - | - | - | - | - | - | - | - | - | - | - | - | - |
| 3 | + | + | + | + | + | + | + | + | + | + | + | + | + | + | + | + | + | + |
| 5 | + | + | + | + | + | + | + | + | + | + | - | - | - | - | - | - | - | - |
| 7 | + | + | + | + | - | - | - | - | - | - | + | + | + | + | + | + | - | - |
| 9 | + | - | - | - | + | + | + | - | - | - | + | + | + | - | - | - | + | + |
| 11 | - | + | - | - | + | - | - | + | + | - | + | - | - | + | + | - | + | + |
| 13 | - | - | + | - | - | + | - | + | - | + | - | + | - | + | - | + | + | - |
| 15 | - | - | - | + | - | - | + | - | + | + | - | - | + | - | + | + | - | + |

Table 3c: Alternatives No. 36 to 53 for e = 16

| 16n + o | Number of alternatives for $^jm_i$ | | | | | | | | | | | | | | | | |
|---|---|---|---|---|---|---|---|---|---|---|---|---|---|---|---|---|---|
| o | 54 | 55 | 56 | 57 | 58 | 59 | 60 | 61 | 62 | 63 | 64 | 65 | 66 | 67 | 68 | 69 | 70 |
| 1 | - | - | - | - | - | - | - | - | - | - | - | - | - | - | - | - | - |
| 3 | + | + | - | - | - | - | - | - | - | - | - | - | - | - | - | - | - |
| 5 | - | - | + | + | + | + | + | + | + | + | + | + | - | - | - | - | - |
| 7 | - | - | + | + | + | + | + | + | - | - | - | - | + | + | + | + | - |
| 9 | + | - | + | + | + | - | - | - | + | + | + | - | + | + | + | - | + |
| 11 | - | + | + | - | - | + | + | - | + | + | - | + | + | + | - | + | + |
| 13 | + | + | - | + | - | + | - | + | + | - | + | + | + | - | + | + | + |
| 15 | + | + | - | - | + | - | + | + | - | + | + | + | - | + | + | + | + |

Table 3d: Alternatives No. 54 to 70 for e = 16

We get an idea for the very fast growth of alternative $^jm_i$ if we list the first three powers of 2 for number e:

| a in e = $2^a$ | Max j in $^jm_i$ of given e |
|---|---|
| 2 | 2 |
| 3 | 6 |
| 4 | 70 |
| 5 | 12870 |

Table 4: Very fast growth for max j

Very interesting are some relations for maximum j:

$$1^2 + 1^2 = \underline{2}$$

$$1^2 + \underline{2}^2 + 1^2 = \underline{6}$$

$$1^2 + 4^2 + \underline{6}^2 + 4^2 + 1^2 = \underline{70} \tag{11}$$

$$1^2 + 8^2 + 28^2 + 56^2 + \underline{70}^2 + 56^2 + 28^2 + 8^2 + 1^2 = 12870$$

The maximum j is the base of the next middle square number:

$$1 + 1 = \underline{2}$$

$$1 + \underline{(1 + 1)}^2 + 1 = \underline{6}$$

$$1 + 4^2 + \underline{(1 + 2^2 + 1)}^2 + 4^2 + 1 = \underline{70}$$

$$1 + 8^2 + 28^2 + 56^2 + \underline{(1 + 4^2 + 6^2 + 4^2 + 1)}^2 + 56^2 + 28^2 + 8^2 + 1 = 12870$$

Another interesting relation occurs if we add in (11) only the base numbers of the square numbers. We find:

$$1 + 1 = 2 = 2^1$$
$$1 + 2 + 1 = 4 = 2^2$$
$$1 + 4 + 6 + 4 + 1 = 16 = 2^4$$
$$1 + 8 + 28 + 56 + 70 + 56 + 28 + 8 + 1 = 256 = 2^8$$

(12)

All these sums are powers of two. The finding of all combinations for changing the leading signs to build quasi alternating prime number series causes those relations because of the both possible signs. The change of a leading sign illustrates the situation about a third case hidden behind the scene.

To finish this short discussion we have a table linking up the special e numbers in (1) from table 4 and the base number sums of (12)

| $e = 2^a$ | Base sums of square numbers calculating max. j in (12) |
|---|---|
| $4 = 2^2$ | $2 = 2^1$ |
| $8 = 2^3$ | $4 = 2^2$ |
| $16 = 2^4$ | $16 = 2^4$ |
| $32 = 2^5$ | $256 = 2^8$ |

Table 5: Growth of base sums

## Appendix B

We do not have only "true" or "false", minus or plus. In last case of minus or plus we know about the third stage to separate both alternatives by zero. The virtual subsets we are talking about are more than a "zero". In quantum physics we discuss a superposition of states. Like particle spin, or in general wave and particle dualism of physical entities and so on. We need an action to decide between the spin direction, the wave and particle behaviour. And that is the observation of physical entities: The measurement! This measurement action in quantum physics we can take as an analogy in mathematics: For observation – measurement – we need to construct the useful instruments for getting experimental results. In mathematics the same action is to construct mathematic objects.

A good example is the discussed possible boundary β used in (5). If a β exists but we can only show formal that there a β must exist. Without construction his explicit value calculation we do have a duality state. On the one hand we know threw formal arguments that there a β must exist but on the other hand we do not know how long we have to search for zero in a given quasi alternate prime number series. The boundary is not sure. That means that β exists potentially and it exists not as a search boundary we need. It is like a superposition of two states. Only a constructed and calculated value of β gives the concrete state of β. But the situation is deeper than that. If we know that there is no β but a series like in (9) then we can not be sure of the existence of a zero. So we have an uncertainty behind an uncertainty: Neither to determinate the value nor the pure existence of a mathematical object. It will be interesting to have also a look in quantum physics to both levels of uncertainties.